\documentclass{amsart}[12pt]
\usepackage{amsmath,amssymb,amsthm,amscd}
\newtheorem{thm}{Theorem}

\newtheorem*{conjec}{Conjecture}

\newtheorem*{ack}{Acknowledgements}

\newcommand{\Z}{\mathbb{Z}}
\newcommand{\legen}[2]{\genfrac{(}{)}{}{}{#1}{#2}}

\begin{document}

\title{A 60,000 digit prime number of the form $x^{2} + x + 41$}
\author{Justin DeBenedetto}
\address{Department of Mathematics, Wake Forest University, Winston-Salem, NC 27109}
\email{debejd0@wfu.edu}

\author{Jeremy Rouse}
\address{Department of Mathematics, Wake Forest University,
  Winston-Salem, NC 27109}
\email{rouseja@wfu.edu}
\thanks{The first author was supported by a Wake Forest Research Fellowship;
the second author was supported by NSF grant DMS-0901090}

\begin{abstract}
Motivated by Euler's observation that the polynomial $x^{2} + x + 41$
takes on prime values for $0 \leq x \leq 39$, we search for large values
of $x$ for which $N = x^{2} + x + 41$ is prime. To apply classical primality
proving results based on the factorization of $N-1$, we choose $x$ to
have the form $g(y)$, chosen so that $g(y)^{2} + g(y) + 40$ is reducible.
Our main result is an explicit, 60,000 digit prime number of the form
$x^{2} + x + 41$.
\end{abstract}

\maketitle

\section{Introduction and Statement of Results}

In 1772, Euler wrote to Johann Bernoulli and mentioned his observation that
$f(x) = x^{2} + x + 41$ takes on prime values for $0 \leq x \leq 39$. Even
after this point, $f(x)$ continues to take on a high frequency of prime
values. For instance, among the numbers $f(1)$, $f(2)$, $f(3)$, $\ldots$,
$f(10^6)$, $261080$ of them are prime. This is more than three times the
number of primes in the sequence $1$, $2$, $3$, $\ldots$, $10^{6}$.

In 1913, Rabinowitsch \cite{Rab} proved that
$n^{2} + n + A$ is prime for $0 \leq n < A-1$ if and only if the ring
\[
  \Z\left[\frac{1 + i \sqrt{D}}{2}\right]
\]
is a principal ideal domain, where $D = 4A - 1$. It is a very
deep result of Baker, Heegner, and Stark that if
$D \equiv 3 \pmod{4}$, then
\[
  \Z\left[\frac{1 + i \sqrt{D}}{2}\right] \text{ is
a principal ideal domain } \iff D = 3, 7, 11, 19, 43, 67 \text{ and } 163.
\]
In \cite{Cox}, pg. 271-274, Cox gives an account of the proof and an overview
of the history of this result.

Euler's polynomial is not unique in taking on a long string of prime values.
Indeed, the polynomial $36x^{2} - 810x + 2753$ (discovered by R. Ruby, see
\cite{Ribenboim}, pg. 112) takes on 45 distinct consecutive prime values
(in absolute value) for $0 \leq x \leq 44$. However, for polynomials
of the form $x^{2} + x + A$, Euler's polynomial still holds the record.

In 1923, Hardy and Littlewood \cite{HL} stated a number of precise conjectures
about the distribution of primes satisfying various additional conditions.
Their prime $k$-tuples conjecture implies that for any positive integer $m$,
there is a number $A$ so that
\[
  x^{2} + x + A \text{ is prime for } 0 \leq x \leq m.
\]
In other words, with a large enough choice of $A$, Euler's polynomial
can be beaten. In addition, they stated a precise conjecture about
how frequent prime values of a fixed quadratic polynomial $f(x)$ are.
\begin{conjec}[See Conjecture F in \cite{Ribenboim}, pg. 190]
Let $a, b, c \in \Z$ with $a > 0$, $\gcd(a,b,c) = 1$, $b^{2} - 4ac$ not
a square where $a+b$ and $c$ are not both even. Let $f(x) = ax^{2} + bx + c$
and let $\pi_{f}(x) = \# \{ p \leq x : p \text{ is prime and } p = f(n) \text{
for some } n \in \Z \}$. Then,
\[
  \pi_{f}(x) \sim \frac{\epsilon C}{\sqrt{a}} \frac{\sqrt{x}}{\log(x)}
  \prod_{\substack{p > 2 \\ p | \gcd(a,b)}} \frac{p}{p-1}.
\]
Here
\[
  \epsilon = \begin{cases} 1 & \text{ if } a+b \text{ is odd }\\
  2 & \text{ if } a+b \text{ is even, }
\end{cases}
\]
and
\[
  C = \prod_{\substack{ p > 2 \\ p \nmid a }} \left(1 - \frac{\legen{b^{2} - 4ac}{p}}{p-1}\right),
\]
and $\legen{b^{2} - 4ac}{p}$ denotes the Legendre symbol.
\end{conjec}

For $f(x) = x^{2} + x + 41$, the conjecture predicts that $\pi_{f}(x)
\sim (6.6395464...) \cdot \frac{\sqrt{x}}{\log(x)}$, and the large
value of the constant $C$ arises
because the values of $x^{2} + x + 41$ are never multiples of primes
$p < 41$. See the papers \cite{FungWilliams} and
\cite{JacobsonWilliams} for computations of larger values of $A$ for
which the corresponding value of the constant $C$ is large.

The goal of the present paper is to provide some verification of the conjecture
of Hardy and Littlewood by searching for large prime values of $x^{2} + x + 41$.
To state our main result, recall that $n\#$ is the product of primes
less than or equal to $n$.

\begin{thm}
Let $f(x) = x^{2} + x + 41$ and $g(x) = 40x^{3} + 41x^{2} + 42x + 1$.
If we set
\[
  x = \frac{310927391 \cdot 23143\#}{43},
\]
then $f(g(x))$ is a $60,000$ digit prime number.
\end{thm}

The fastest current methods for proving primality of large numbers $N$
are based on knowing partial prime factorizations of $N-1$ or $N+1$. The
best general primality proving method not based on factorizations is
the elliptic curve primality proving method (ECPP), and the current record for
a general number is 26,642 digits. This number was proven prime by
Fran\c cois Morain in 2011. (Note that in \cite{ASSW}, primality of
numbers in a very particular sequence is proven using ECPP. Some of these
numbers have more than 100,000 digits, but this method does not apply in
general.)

Prime numbers of the form $x^{2} + 1$ have received special attention,
and the largest known such prime is (as of this writing)
$75898^{524288} + 1$, with 2,558,647 digits. It is straightforward to
find large primes of this type since for a number of the
form $x^{2} + 1$ we can factor $N-1$ as long as we know the prime factorization
of $x$.

This is not the case for $f(x) = x^{2} + x + 41$. Our approach to finding large
primes of this type is to find polynomials $g(x)$ so that $f(g(x)) - 1$
is reducible. A computer search revealed the choice
$g(x) = 40x^{3} + 41x^{2} + 42x + 1$ for which
\[
  f(g(x)) - 1 = (40x^{2} + x + 1) (40x^{4} + 81x^{3} + 123x^{2} + 84x + 42).
\]
The Brillhart-Lehmer-Selfridge theorem (see \cite{BLS}, Theorem~5 or
Section~2) allows one to prove the primality of $N$ provided we know
the complete factorization of a factor $F$ of $N-1$ with $F$ larger
than about $N^{1/3}$. Our goal then is to find a choice of $x$ (with
known prime factorization) for which $40x^{2} + x + 1$ is prime (proven again
using the Brillhart-Lehmer-Selfridge theorem), and for which $f(g(x))$ is
prime. This simultaneous primality requirement significantly increases
the number of candidate values of $x$ we must search, and makes our
result comparable in difficulty to finding a large twin prime
or Sophie Germain prime.

An outline of the paper is as follows. In the second section, we give
appropriate background, and in the third section we describe our computations
and the verification of the primality of our 60,000 digit number of the
form $x^{2} + x + 41$.

\begin{ack}
We used PARI/GP \cite{PARI} for sieving computations, OpenPFGW
for primality testing, and the Wake Forest DEAC cluster for primality
testing computations. We would like to thank David Chin for compiling
OpenPFGW for us on the DEAC cluster.
\end{ack}

\section{Strategy}

We start by stating the Brillhart-Lehmer-Selfridge theorem.

\begin{thm}[A special case of Theorem~5 of \cite{BLS}]
Suppose that $N > 1$ is odd and write $N - 1 = FR$ where $F$ is even and the prime factorization of $F$ is known.  Suppose also that
\begin{enumerate}
\item $F > (\frac{N}{2})^{1/3}$,
\item For each prime $p_{i}$ dividing $F$, there is an integer $a_i$ so that ${a_i}^{N - 1} \equiv 1 \pmod{N}$ and gcd$({a_i}^{\frac{N - 1}{{p_i}}} - 1, N) = 1$,
\item If we write $R = 2Fq + r$, where $1 \leq r < 2F$, then either $q = 0$ or $r^2 - 8q$ is not a perfect square,
\end{enumerate}
then $N$ is prime.
\end{thm}

In order to take advantage of the Brillhart-Lehmer-Selfridge Theorem for proving primality, we used the following equations\\
\begin{align*}
f(g(x)) - 1 &= h(x)i(x)\\
h(x) &= 40x^2 + x + 1\\
i(x) &= 40x^4 + 81x^3 + 123x^2 + 84x + 42,\\
\end{align*}
with $f(x)$ and $g(x)$ defined above.  Thus for any given choice of $x$, with $x$ even, we have $N = f(g(x))$ with $2h(x) = F$ and $\frac{i(x)}{2} = R$.

Next, we estimate how many numbers we will have to test.  According to the Prime Number Theorem, the density of primes close to an integer $N$ is approximately equal to $\frac{1}{\ln(N)}$.  Since we were looking for a 20,000 digit number and the corresponding 60,000 digit number to be simultaneously prime, if we assume the same density of primes within the values of Euler's Polynomial as within the set of all integers, we must multiply the probability of finding a 20,000 digit prime with the probability of finding a 60,000 digit prime.  Thus our expected probability for a given pair of values (corresponding to $h(x)$ and $f(g(x))$ from above) is approximately
\begin{center}
$(\frac{1}{\ln(10^{20000})})(\frac{1}{\ln(10^{60000})}) \approx \frac{1}{6,362,314,060}$.
\end{center}
If we were to test 6,362,314,060 numbers, our chance of finding at least one prime pair would be
\begin{center}
$1 - (\frac{6,362,314,059}{6,362,314,060})^{6,362,314,060}$.
\end{center}
We know that $1 - ((N-1)/N)^N \approx 1 - \frac{1}{e}$, thus
\begin{center}
$1 - (\frac{6,362,314,059}{6,362,314,060})^{6,362,314,060}\approx 1 - \frac{1}{e}\approx 63.2\%$.\\
\end{center}
We chose to broaden our search in order to have a higher theoretical probability of success.  We wanted to be roughly 95\% confident that our search will yield success, thus we tripled our amount of numbers to check.
\begin{center}
$1 - (\frac{6,362,314,059}{6,362,314,060})^{3 \cdot 6,362,314,060} \approx 1 - (\frac{1}{e})^3 \approx 95.0\%$.
\end{center}
Since the equation we were working with was
\begin{center}
$f(g(x)) = 1600x^6 + 3280x^5 + 5041x^4 + 3564x^3 + 1887x^2 + 126x + 43,$
\end{center}
we chose to use a primorial divided by $43$ as our $x$.  By working with an integer multiple of a primorial as our value for $x$, $x = k\cdot\frac{n\#}{43}$, we know that $f(g(x))$ will not be divisible by any prime less than $n$.  This is because we know that if $x \equiv 0 \pmod{p}$,
then $h(x) \equiv 1 \pmod{p}$ and $f(g(x)) \equiv 43 \pmod{p}$.  Since $f(g(x)) \equiv 43 \pmod{p}$ if $x \equiv 0 \pmod{p}$, we do not want $x$ to be a multiple of $43$.

For each potential prime divisor eliminated, the number of potential primes decreases by $\frac{1}{p}$ where $p$ is the divisor eliminated.  This is due to the fact that the density of numbers $n$ divisible by $p$ is $\frac{1}{p}$.  Thus the number of numbers we should check should be
\begin{center}
$(\prod_{\text{p prime}<23,143} \frac{p-1}{p})^2 (3 \cdot 6,362,314,060) \approx 59,481,223$,
\end{center}
due to our use of $\frac{23,143\#}{43}$ as a factor of $x$.

Next, we chose to further reduce our search by sieving the test numbers.  We eliminated all numbers that were divisible by primes under $5\cdot10^{10}$.  The number of numbers left after sieving up to $5\cdot10^{10}$ is approximately
\[
  3 \cdot (6,362,314,060) \cdot
  \prod_{\substack{p \text{ prime}\\ p \leq 5\cdot10^{10}}} \left(1 - \frac{1}{p}\right)^{2}.
\]
Mertens's theorem states that
\[
  \prod_{\substack{p \text{ prime}\\ p \leq x}} \left(1 - \frac{1}{p}\right)
  \sim \frac{1}{e^{\gamma} \ln(x)}.
\]
Using this approximation, we estimate that 9,914,204 numbers would remain after sieving up to $5\cdot10^{10}$.

Finally, we estimate how much CPU time we will need to use.  On our computers, primality tests took approximately 14 seconds for 20,000 digit numbers, and 123 seconds for 60,000 digit numbers. This left us with the following preliminary CPU time estimations\\
Without utilizing primorials:\\ $3 \cdot 6,362,314,060 \cdot 14 \text{ seconds} + 3 \ln(10^{60000}) \cdot 123 \text{ seconds} \approx 8,475 \text{ years}$.\\
Using Primorials Pre-Sieve:\\ $59,481,223 \cdot 14 \text{ seconds} + 59,481,223 / \ln(10^{20000}) \cdot 123 \text{ seconds} \approx 26 \text{ years}$.\\
After Sieving:\\ $9,914,204 \cdot 14 \text{ seconds} + 9,914,204 / \ln(10^{20000}) \cdot 123 \text{ seconds} \approx 4 \text{ years}$.\\

\section{Computations}

We begin by describing how the sieving was done.  In $f(g(x))$ we had $x = k(\frac{23,143\#}{43})$, thus the sieving was designed to eliminate values of $k$ for which $f(g(x))$ was a multiple of a prime $p$.  To do this, we looked for roots of the polynomials, $f(x)$ and $f(g(x))$, in $\mathbb{F}_p$.  We implemented the sieving through a program we wrote to run through Pari/GP.  For each prime, $p$, in our testing range (below $5 \cdot 10^{10}$) we took the following steps.  First, we factored $40x^2+x+1$ mod $p$ and computed its roots in $\mathbb{F}_p$.  Then, we eliminated the choices of $k$ for which the corresponding $x$ value is a root.  Once we had completed this process, we then repeated this process on $1600x^6+3280x^5+5041x^4+3564x^3+1887x^2+126x+43$, once again eliminating roots in $\mathbb{F}_p$.

Next, we examine our sieving results.  We sieved up to $5 \cdot 10^{10}$ and brought our total number of numbers down from 60,000,000 to 9,946,272.  This is just over 30,000 more than we had estimated would be left after sieving.  This process took approximately five days running on a single computer.

Next, we describe how the pseudo-primality testing was done.  We ran Fermat pseudo primality tests to find probable primes base 3 on OpenPFGW.  OpenPFGW uses the fast Fourier transform (FFT) method for fast multiplication.  We ran our tests in groups of around 900 numbers on the DEAC cluster, which has approximately 1,200 nodes with processors ranging from 2.4 GHz to 3.0 GHz.  We automated job submission to the cluster by checking how many jobs were currently running and how many nodes were in use and making appropriate choices for submitting more jobs based on that information.  When a pseudo-prime was found, that number was appended to a file, giving us a single list of all 20,000 digit pseudo-primes found.  While the tests were running, we periodically checked the corresponding 60,000 digit numbers for pseudo-primality base 3.

Finally, we examine our computation results.  We tested a bit more than 3,000 of the resulting 60,000 digit numbers and found one pseudo-prime.  We stopped running tests once we successfully confirmed through OpenPFGW that the 60,000 digit pseudo-prime was prime.  We had found 3,521 20,000 digit pseudo-primes, and the 2,813th one corresponded to our 60,000 digit prime.  Of the 9,946,271 numbers post-sieving, the 20,000 digit number yielding the 60,000 digit prime was the 3,106,282nd.  The total amount of CPU time that was used was: $5$ days for sieving, $14 \cdot (3,300,000)$ seconds for 20,000 digit pseudo-primality tests, and $123 \cdot (3,000)$ seconds for 60,000 digit pseudo-primality tests, totalling approximately $544$ days.

\bibliographystyle{plain}
\bibliography{refs}

\end{document}